\documentclass[12pt]{amsart}

\usepackage{amssymb,amsfonts} 	%AMS Schriften & Symbole
\usepackage{amsmath}            %AMS Mathematische Befehle wie align...
\usepackage{amsthm}
\usepackage{fullpage}
\usepackage{pstricks}
\usepackage{booktabs}

\newtheorem{thm}{Theorem}
\newtheorem{lem}{Lemma}
\newtheorem{prop}{Proposition}
\newtheorem{cor}{Corollary}
\newtheorem{conj}[thm]{Conjecture}

\theoremstyle{remark} % 'style changed again'
\newtheorem{rem}{Remark}
\theoremstyle{definition}

\def\modd#1 #2{#1\ ({\rm mod}\ #2)}

\renewcommand{\le}{\leqslant}
\renewcommand{\ge}{\geqslant}
\renewcommand{\Xi}{\varXi}

\newcommand{\eps}{\varepsilon}

\newcommand{\N}{\mathbb{N}}

\begin{document}

\title{Thue-Morse at Multiples of an Integer}

\author{Johannes F. Morgenbesser}
\address{Institut f\"ur Diskrete Mathematik und Geometrie, Technische Universit\"at Wien, Wiedner Hauptstra\ss e 8--10, A--1040 Wien, Austria,
\newline \indent
Institut de Math\'ematiques de Luminy, Universit\'e de la M\'editerran\'ee, 13288 Marseille Cedex 9, France,}
\email{johannes.morgenbesser@tuwien.ac.at}
\thanks{J. Morgenbesser was supported by the Austrian Science Foundation FWF, grant S9604, that is part of the National Research Network ``Analytic Combinatorics and Probabilistic Number Theory''.}

\author{Jeffrey Shallit}
\address{School of Computer Science, University of Waterloo, Waterloo, ON  N2L 3G1, Canada,}
\email{shallit@cs.uwaterloo.ca}

\author{Thomas Stoll}
\address{Institut de Math\'ematiques de Luminy, Universit\'e de la M\'editerran\'ee, 13288 Marseille Cedex 9, France,}
\email{stoll@iml.univ-mrs.fr}

\keywords{Thue-Morse sequence, sum of digits, congruences, arithmetic progressions}
\subjclass[2010]{11N25, 11A63, 68R15}
\maketitle

\begin{abstract}
Let ${\bf t} = (t_n)_{n\geq 0}$ be the classical Thue-Morse sequence
defined by $t_n = s_2(n) \bmod 2$, where $s_2$ is the sum of the bits
in the binary representation
of $n$. It is well known that for any integer $k \geq 1$ the frequency
of the letter ``1'' in the subsequence $t_0, t_k, t_{2k}, \ldots$ is
asymptotically $1/2$. Here we prove that for any $k$ there is a $n\leq
k+4$ such that $t_{kn}=1$. Moreover, we show that $n$ can be chosen to
have Hamming weight $\leq 3$. This is best in a twofold sense. First,
there are infinitely many $k$ such that $t_{kn}=1$ implies that $n$
has Hamming weight $\geq 3$. Second, we characterize all $k$ where
the minimal $n$ equals $k$, $k+1$, $k+2$, $k+3$, or $k+4$. Finally, we
present some results and conjectures for the generalized problem, where
$s_2$ is replaced by $s_b$ for an arbitrary base $b\geq 2$.

\end{abstract}

%%%%%%%%%%%%%%%%%%%%%%%%%%%%%%%%%%%%%%%%%%%%%%%%%%%%%%%%%%%%%%%%%%%%%%%%%%%%%%%%%%%%%%%%%%%%%%%%%%%%%%%%%%%%%%%%%%%%%%%%%%%%%%%%%%%%%%
%%%%%%%%%%%%%%%%%%%%%%%%%%%%%%%%%%%%%%%%%%%%%%%%%%%%%%%%%%%%%%%%%%%%%%%%%%%%%%%%%%%%%%%%%%%%%%%%%%%%%%%%%%%%%%%%%%%%%%%%%%%%%%%%%%%%%%

\section{Introduction}

Let $s_b (n)$ denote the sum of the digits of $n$ when expressed in
base $b$, and let $$t_n = s_2 (n) \bmod 2,\qquad n\geq 0,$$ be the
Thue-Morse sequence $\bf t$. In an e-mail message dated June 7 2010,
Jorge Buescu of the Universidade de Lisboa observed that the Thue-Morse
sequence can be regarded as a $2$-coloring of the integers, and
therefore, by van der Waerden's theorem, must contain arbitrarily long
monochromatic arithmetic progressions.\footnote{We do not need the
power of van der Waerden's theorem to prove this. For example, as we
will see later, if $k = 2^r - 1$ for some $r \geq 1$, then $s_2 (kn) = r$
for $1 \leq n \leq k$.} (By a {\it monochromatic arithmetic
progression} we mean a series of indices $i, i+j, i+2j, \ldots,
i+(n-1)j$ such that $t_i = t_{i+j} = \cdots = t_{i+(n-1)j}$.) He then
asked, is it true that $\bf t$ has no infinite monochromatic arithmetic
progressions?

The answer is yes:  $\bf t$ has no infinite monochromatic arithmetic
progressions. This is a consequence of the following result of Gelfond
\cite{Gelfond:1968}, which says that the values of $s_b (n)$ are
equally distributed in residue classes, even if the residue class of
$n$ is fixed.  (A weaker result, applicable in the case of the
Thue-Morse sequence, had previously been given by Fine~\cite{Fine:1965}.)

\begin{thm}
Let $b, r, m$ be positive integers with $\gcd(b-1, r) = 1$, and let $a, c$ be any integers.  Then the number of integers $n \leq x$ congruent to $a$ mod $k$ such that
$s_b (n) \equiv \modd{c} {r}$ is equal to $\frac{x}{kr} + O(x^\lambda)$ for some $\lambda < 1$ that does not depend on $x, k, a$, or $c$.
\end{thm}

Gelfond's theorem, however, concerns the {\it average} distribution of the values of $s_b (n)$ in residue classes.  It suggests the following question:  how large can the smallest $n$ be that is
congruent to $a$ mod $k$ and satisfies $s_b (n) \equiv \modd{c} {r}$?

In this paper we answer the question for the case $a=0$, $k$ arbitrary,
$c=1$, $b=r=2$. In other words, we find a bound on the number of terms
in a fixed arithmetic progression of the Thue-Morse sequence we have to
look at in order to see a ``1''. We include some weaker results for
arbitrary $b$ and give some conjectures.

\begin{rem}
Jean-Paul Allouche notes that Buescu's original question can also be answered by appealing to a lemma in his paper \cite[p.\ 284]{Allouche:1982}.  His lemma states that if $a, b, c$ are integers with $b - c > a$, then $t_{an+b} - t_{an+c}$ cannot be constant for large $n$. If $t_{An+B}$ were constant for some integers $A, B$ then it would have the same value when replacing $n$ by $n+2$. Thus $t_{An + 2A + B} - t_{An+B}$ would be constant and equal to $0$, but $2A+B-B = 2A > A$ and we are done.
\end{rem}

\begin{rem}
  Dartyge, Luca and St{\u{a}}nic{\u{a}}~\cite{DartygeLucaStanica:2009} recently investigated another problem on the {\it pointwise} behavior of $s_b$ on integer multiples, namely, to bound the smallest nontrivial $n$ that is congruent to $0$ mod $k$ and satisfies $s_b(n) = s_b(k)$. For other distributional properties of $s_b$ on integer multiples we refer the interested reader to the bibliographic list in~\cite{DartygeLucaStanica:2009}.
\end{rem}

\medskip

To begin with, for $k\geq 1$ we write $$\mathcal{N}_k=\{n: t_{kn} =1\},\qquad
f(k) = \min\{ n : n\in \mathcal{N}_k \}.
$$
The first few values of $(f(k))_{k\ge1}$ are given by
\begin{align}\label{eq:5}
1,\;1,\;  7,\; 1,\; 5,\; 7, \; 1, \; 1, \; 9,\; 5,\; 1, \;7,\;1,\; 1,\; 19,\; 1,\;  17,\; 9,\; 1,\;5\ldots
\end{align}
From Gelfond's theorem we get that $f(k)<\infty$ for all $k$.
Indeed, a simple observation
shows that $f(k)=O(k)$. To see this, consider $n=2^{2r+1}-1\geq k$.
Then
$$s_2(kn)=s_2(k 2^{2r+1} - k)= s_2(k-1)+2r+1-s_2(k-1) \equiv \modd{1} {2},$$
so that $f(k)< 4k$.

\medskip

This function is of interest because of some old work of Newman \cite{Newman:1969}. Leo Moser observed that the first few multiples of $3$ all have an even number of digits in their base-$2$ expansion. In our notation, this means $f(3) = 7$.  In fact, Newman showed that among the first $x$ multiples of $3$, there is always a small preponderance of those with even parity.

Our main result shows that $k=3$ is the first of an infinite class of integers that maximize $f(k)-k$.

\begin{thm}\label{thm:1}
For all $k\geq 1$ we have
\begin{align}\label{eq:1}
f(k)\le k+4.
\end{align}
Moreover, we have
\begin{enumerate}
\item $f(k)=k+4$ if and only if $k=2^{2r}-1$ for some $r\ge1$,
\item There are no $k$ with $f(k)=k+3$ or $f(k)=k+2$.
\item $f(k)=k+1$ if and only if $k=6$.
\item $f(k)=k$ if and only if $k=1$ or $k=2^r + 1$ for some $r\ge2$.
\end{enumerate}
\end{thm}

In fact, we can always find such a small $n$ with Hamming weight at most 3.

\begin{cor}
  For all $k$ there is an $n\in \mathcal{N}_k$ with $n\leq k+4$ and $s_2(n)\leq 3$.
\end{cor}

This is optimal in the sense that there are infinitely many $k$ such that all $n\in \mathcal{N}_k$ satisfy $s_2(n)\geq 3$.
To see this, consider $k=3\cdot 2^r+3$ for $r\geq 4$. Then $s_2(k\cdot 1)=4\equiv \modd{0} {2}$ and $s_2(k(2^j+1))\in \{4,6,8\}$ for $j\geq 1$.

\medskip

The paper is structured as follows. We introduce some useful notation in Section~\ref{section:notation} which allows us to perform addition in the binary expansion
of integers in a well-arranged way. In Section~\ref{section:idea} we shortly outline the idea of the proof of our main result.
In Section~\ref{section:aux} we state some auxiliary results which are based on a detailed investigation of various cases. Section~\ref{section:theorem}
is devoted to the proof of Theorem~\ref{thm:1}. We conclude with some results in the general case, where the condition $s_2(kn)\equiv \modd{1} {2}$
is changed to $s_b(kn)\equiv \modd{c} {r}$ (Section~\ref{section:general}).

%%%%%%%%%%%%%%%%%%%%%%%%%%%%%%%%%%%%%%%%%%%%%%%%%%%%%%%%%%%%%%%%%%%%%%%%%%%%%%%%%%%%%%%%%%%%%%%%%%%%%%%%%%%%%%%%%%%%%%%%%%%%%%%%%%%%%%
%%%%%%%%%%%%%%%%%%%%%%%%%%%%%%%%%%%%%%%%%%%%%%%%%%%%%%%%%%%%%%%%%%%%%%%%%%%%%%%%%%%%%%%%%%%%%%%%%%%%%%%%%%%%%%%%%%%%%%%%%%%%%%%%%%%%%%

\section{Notation}\label{section:notation}

In this section we introduce some notation.
If
\[
k = \eps_{\ell-1}(k) 2^{\ell-1} + \eps_{\ell-2}(k) 2^{\ell-2} + \cdots + \eps_0(k)
\]
is the canonical base-$2$ representation of $k$, satisfying
$\eps_{j}\in\{0,1\}$ for all $0\le j <\ell$ and $\eps_{\ell-1}(k)\ne 0$,
then we let $(k)_2$ denote
the binary word $$ \eps_{\ell-1}(k)  \eps_{\ell-2}(k)\cdots \eps_0(k).$$ Additionally, for each $k\in\N$  we let $\ell(k)$ denote the length of
$(k)_2$; for $k \geq 1$ this is $\ell(k) = \lfloor \log_2 k \rfloor +1$.
If $w_1$ and $w_2$ are two binary words, then $w_1w_2$ denotes the binary word obtained by concatenation. The symbol $a^n$, $n\ge1$, $a\in\{0,1\}$
is an abbreviation for the word
\[
\overbrace{aa\cdots a}^{n},
\]
and $a^0$ is equal to the empty word. For $a\in\{\,0\,,1\}$, we use the notation $\bar{a} = 1 - a$.
We define the function $s$ for all binary words $w= \eps_{j-1}\cdots\eps_0$ by
\[
s(w) = \# \{ 0\le i < j : \eps_i=1\},
\]
and in particular, we have $s((k)_2) = s_2(k) \equiv \modd{t_k} {2}$. If $1\le j\le \ell(k)$, we set
\[
L_j(k) = \eps_{j-1} \cdots \eps_0,
\]
the $j$ least significant bits of $k$ in base $2$,
and
\[
U_j(k) = \eps_{\ell(k)-1} \cdots \eps_{\ell(k)-j},
\]
the $j$ most significant bits of $k$.
For example, if $k=119759$, then we have $(k)_2 = 1\, 1 \,1\, 0\, 1 \,0\, 0\, 1\, 1\, 1\, 1\, 0\, 0\, 1\, 1\, 1\, 1$, $\ell(k)=17$ and
\[
\overbrace{1\; 1\; 1\; 0\; 1\; 0\; 0\; 1\; 1\;
\makebox[0pt][l]{$\displaystyle{\underbrace{\phantom{1\; 1\; 0\;} 0 \;1\; 1\; 1\; 1}_{L_8(k)}}$}
1\; 1\; 0}^{U_{12}(k)}\phantom{\;0 \;1\; 1\; 1\; 1}.
\]
Written in short form, this means that $(k)_2 = 1^3\, 0\, 1\, 0^2\, 1^4\, 0^2\, 1^4$,
\[
U_{12}(k) = 1^3\, 0\, 1\, 0^2\, 1^4\, 0 \quad \mbox{ and }\quad L_8(k) = 1^2 \, 0^2 \, 1^4.
\]
In what follows, we use the convention that if we are talking about $L_j(k)$ or $U_j(k)$, we assume that $\ell(k)\ge j$.
Note that for all $k\in\N$ and $j< \ell(k)$ we have
\[
s(L_{\ell(k)-j}(k)) \equiv \modd{s(U_j(k)) + t_k} {2} .
\]
Furthermore, this function also satisfies
\[
s(w_1 w_2) = s(w_1) + s(w_2)
\]
for two binary words $w_1$ and $w_2$.

%%%%%%%%%%%%%%%%%%%%%%%%%%%%%%%%%%%%%%%%%%%%%%%%%%%%%%%%%%%%%%%%%%%%%%%%%%%%%%%%%%%%%%%%%%%%%%%%%%%%%%%%%%%%%%%%%%%%%%%%%%%%%%%%%%%%%%
%%%%%%%%%%%%%%%%%%%%%%%%%%%%%%%%%%%%%%%%%%%%%%%%%%%%%%%%%%%%%%%%%%%%%%%%%%%%%%%%%%%%%%%%%%%%%%%%%%%%%%%%%%%%%%%%%%%%%%%%%%%%%%%%%%%%%%

\section{Idea of proof}
\label{section:idea}

It is relatively easy to show that $f(k)=k+4$ if $k=2^{2r}-1$ for some $r\ge1$ and $f(k)=k$ if $k=2^{n}+1$ for some $n\ge2$ (see the proofs of Theorem~\ref{thm:1} and Lemma~\ref{lem:1}).
Moreover, since $f(k)=f(2k)$ for all $k\ge1$,
in order to prove Theorem~\ref{thm:1},
it suffices to show that
$f(k)< k$ for all odd integers $k$ other than those stated above.
Thus, we assume in
Section~\ref{section:aux} that $k$ is an odd integer.

We use two different ideas in order to succeed,
depending on the base-$2$ representation of $k$.
We show for a large set of integers $k$ that there exists an integer $n<k$ with Hamming weight two such that $t_{kn}=1$.
To be more precise, we find for such integers $k$ a positive integer $a<\ell(k)$
such that $t_{kn}=1$ with $n= 2^{a}+1\le 2^{\ell(k)-1}+1 <k$.
For the remaining
odd integers $k$ we show that there exist positive odd integers $m<k$ with Hamming weight two and $n<k$ with Hamming weight three such that
\[
t_{kn} \equiv \modd{1+ t_k + t_{km}} {2} .
\]
This implies that $f(k)<k$ since at least one of the three numbers $t_k$, $t_{km}$ and $t_{kn}$ has to be equal to $1$.

%%%%%%%%%%%%%%%%%%%%%%%%%%%%%%%%%%%%%%%%%%%%%%%%%%%%%%%%%%%%%%%%%%%%%%%%%%%%%%%%%%%%%%%%%%%%%%%%%%%%%%%%%%%%%%%%%%%%%%%%%%%%%%%%%%
%%%%%%%%%%%%%%%%%%%%%%%%%%%%%%%%%%%%%%%%%%%%%%%%%%%%%%%%%%%%%%%%%%%%%%%%%%%%%%%%%%%%%%%%%%%%%%%%%%%%%%%%%%%%%%%%%%%%%%%%%%%%%%%%%%

\section{Auxiliary results}
\label{section:aux}

We have to distinguish several cases according to
the beginning and the ending part of the binary expansion of $k$.
\begin{lem}\label{lem:1}
Let $k\in\N$ such that there exists an odd integer $u\ge1$ with $L_{u+1}(k)=01^u$.
Then we have $f(k)\le k$. Furthermore, $f(k)=k$ if and only if $k=2^r+1$ for some $r\ge2$.
\end{lem}

\begin{proof}
Let $\ell=\ell(k)$ and set $n=2^{\ell-1}+1$. In what follows we show that $t_{kn}=1$. We have
\[
(kn)_2 = U_{\ell-(u+1)}(k)\, 1\, 0^u \,L_{\ell-1}(k).
\]
The following figure explains this fact:
\begin{center}
\begin{tabular}{ll}
$\phantom{ .}\cdots \; 0\; 1^{u-1} \;1$ & \\
$\phantom{ .\cdots \; 0\; 1^{u-1}} \;1\; \cdots$ & \\
\cmidrule{1-1}
%\hline
$\phantom{ .}\cdots \; 1 \; 0^{u-1} \; 0\; \cdots\phantom{ .}$& \hspace{-0.12cm}\mbox{\raisebox{+0.42cm}[0pt]{.}}
\end{tabular}
\end{center}

The first line ($\cdots \; 0\; 1^{u-1} \;1$) corresponds to the expansion of $k 2^{\ell-1}$ and the second line ($1 \; \cdots$)
to the expansion of $k$. By ``$\cdots$'' we refer to digits that are not important for our argument.
Since
\[
s(U_{\ell-(u+1)}(k)) \equiv s(L_{u+1}(k)) + t_k \equiv \modd{u +t_k} {2},
\]
and
\[
s(L_{\ell-1}(k)) \equiv s(U_{1}(k)) + t_k \equiv \modd{1 + t_k} {2},
\]
we obtain
\[
t_{kn} \equiv s((kn)_2) \equiv u + t_k  + 1 +1 + t_k \equiv u \equiv
	\modd{1} {2},
\]
which shows that $t_{kn}=1$.
The definition of $\ell=\ell(k)$ implies that $2^{\ell-1}+1\le k$. If $k=2^{\ell-1}+1$,
we have $t_{km} =0$ for all $1\le m <k$. Indeed, if $1\le m < 2^{\ell-1}$, then the
$2$-additivity of the binary sum-of-digits function $s_2$ implies
\[
s_2(km) = s_2(2^{\ell-1}m + m) = s_2(m) + s_2(m).
\]
Thus we have $t_{km} \equiv s_2(km) \equiv \modd{0} {2}$ for all $1\le m < 2^{\ell-1}$.
If $m=2^{\ell-1}$, then we clearly have $t_{km}=0$. This finally proves that $f(k)=k$
if $k=2^{\ell-1}+1$ and $f(k)<k$ if $k$ satisfies the assumptions of Lemma~\ref{lem:1} but $k\ne 2^{\ell-1}+1$.
\end{proof}

%%%%%%%%%%%%%%%%%%%%%%%%%%%%%%%%%%%%%%%%%%%%%%%%%%%%%%%%%%%%%%%%%%%%%%%%%%%%%%%%%%%%%%%%%%%%%%%%%%%%%%%%%%%%%%%%%%%%%%%%%%%%%%%%%%
%%%%%%%%%%%%%%%%%%%%%%%%%%%%%%%%%%%%%%%%%%%%%%%%%%%%%%%%%%%%%%%%%%%%%%%%%%%%%%%%%%%%%%%%%%%%%%%%%%%%%%%%%%%%%%%%%%%%%%%%%%%%%%%%%%

\begin{lem}\label{lem:2}
Let $k\in\N$. If there exists an even integer $u\ge2$ with $L_{u+2}(k)=101^u$,
then we have $f(k)< k$.
\end{lem}

\begin{proof}
Set $\ell=\ell(k)$. First, we show that if there exists a positive integer $r\ne u$ such that $U_{r+1}(k)=1^r 0$, then $f(k)<k$.
If $r<u$, we set $n= 2^{\ell-(r+1)}+1<k$. Then we have
\[
(kn)_2 = U_{\ell-(u+1)}(k)\, 1\, 0^{u-(r+1)}\, 1^{r-1}\, 0\, 1\, L_{\ell-r-1}(k),
\]
as illustrated below:
\begin{center}
\begin{tabular}{ll}
$\phantom{ .}\cdots \; 0\; 1^{u-(r+1)} \;1^{r-1} \; 1 \;1$&\\
$\phantom{ .\cdots \; 0\; 1^{u-(r+1)}} \;1^{r-1} \; 1 \;0\; \cdots$&\\
\cmidrule{1-1}
$\phantom{ .}\cdots \; 1 \; 0^{u-(r+1)}\;1^{r-1} \; 0\;1\; \cdots\phantom{ .}$& \hspace{-0.12cm}\mbox{\raisebox{+0.42cm}[0pt]{.}}
\end{tabular}
\end{center}
Since $s(U_{\ell-(u+1)}(k)) \equiv s(L_{u+1}(k)) + t_k \equiv \modd{u + t_k}
{2}$ and $s(L_{\ell-r-1}(k)) \equiv s(U_{r+1})+t_k \equiv r+ t_k$, we get
\begin{align}\label{eq:3}
t_{kn} \equiv u + t_k + 1 + (r-1) + 1+ r+ t_k \equiv \modd{u+1} {2}.
\end{align}
This shows that $f(k)<k$ if $r<u$ since $u$ is even. If $r>u$, we set
$n=2^{\ell-u} + 2^{\ell-u-1} + 1$. Since $\ell-u<\ell-1$ we have $n<k$. We get
\[
(kn)_2 = U_{\ell+2-(u+2)}(3k)\, 1\, 0 \, 1^{u-2} \, 0^2 L_{\ell-(u+1)}(k),
\]
as illustrated below:
\begin{center}
\begin{tabular}{ll}
$\phantom{ .}\cdots \;1 \; 0\;1\;  1^{u-2} \;1$&\\
$\phantom{ .1\;}\cdots\; 1\; 0 \; 1^{u-2}\; 1 \;1$&\\
$\phantom{ .1\;\cdots\; 1\;}  1\; 1^{u-2}\; 1 \;1\; \cdots$&\\
\cmidrule{1-1}
$\phantom{ .0}\;\cdots\; 1\; 0 \; 1^{u-2}\; 0 \;0\; \cdots$\phantom{ .}& \hspace{-0.12cm}\mbox{\raisebox{+0.42cm}[0pt]{.}}
\end{tabular}
\end{center}
Noting that $3k$ has $\ell+2$ digits, i.e., $\ell(3k)=\ell+2$, we obtain
\begin{align*}
t_{kn} &\equiv s(L_{u+2}(3k)) + t_{3k} + 1+(u-2) + s(U_{u+1}(k)) + t_k\\
 &\equiv s(L_{u+2}(3k)) + t_{3k} +u-1 +u+1 +t_k\\
 &\equiv \modd{s(L_{u+2}(3k)) + t_{3k} +t_k} {2} .
\end{align*}
Since
\begin{align*}
L_{u+2}(3k) = 0^2 1^{u-2}\, 0 \,1,
\end{align*}
we have $t_{kn} \equiv \modd{1 + t_k +t_{3k}} {2}$.
As we have seen in Section~\ref{section:idea}, this implies $f(k)<k$.

For the rest of the proof we assume that $U_{u+1}(k)=1^u 0$ for
some integer $u \geq 2$.
If $(k)_2 = 1^u 0 1^u$, then it is easy to see that $f(k)=3$.
Thus we can assume that there exists a positive integer $v$ such that  $L_{v+u+2}(k) = 0 1^v 0 1^u$. If $v$ is odd, then we set $n=2^{\ell-(u+1)}+1<k$. We get
\[
(kn)_2 = U_{\ell-(v+u+2)}(k)\, 1\, 0^{v+1} \, 1^{u-2} \, 0 \, 1\, L_{\ell-(u+1)}(k),
\]
as illustrated below:
\begin{center}
\begin{tabular}{ll}
$\phantom{ .}\cdots \;0 \; 1^v \;0\;  1^{u-2} \;1\;1$&\\
$\phantom{ . \cdots \;0 \; 1^v} \;1\;  1^{u-2} \;1\;0\;\cdots$&\\
\cmidrule{1-1}
$\phantom{ .}\cdots \;1 \; 0^v \;0\;  1^{u-2} \;0\;1\;\cdots$\phantom{ .}&\hspace{-0.12cm}\mbox{\raisebox{+0.42cm}[0pt]{.}}
\end{tabular}
\end{center}
This implies
\begin{align*}
t_{kn} &\equiv s(L_{v+u+2}(k)) + t_{k} + 1+(u-2) +1 + s(U_{u+1}(k)) + t_k\\
&\equiv u+v + u \equiv \modd{1} {2} .
\end{align*}
If $v$ is even, we have two cases to consider:  $u \ge 4$ and $u = 2$.

If $u\ge4$, we set $n= 2^{\ell-u} + 2^{\ell-u-1}+1<k$. Then we have
\[
(kn)_2 = U_{\ell+2-(v+u+2)}(3k)\,0 \, 1^{v}\, 0 \, 1^{u-3} \, 0 \, 1^2 L_{\ell-(u+1)}(k),
\]
as illustrated below:
\begin{center}
\begin{tabular}{ll}
$\phantom{ .}\cdots \;0 \;1\; 1^{v-1} \;0\;1\;  1^{u-3} \;1\;1$&\\
$\phantom{ .\;0\;}\cdots \; 0\; 1^{v-1} 1\;0\;  1^{u-3} \;1\;1\;1$&\\
$\phantom{ .\;0\;\cdots \; 0\; 1^{v-1} 1\;}1\;  1^{u-3} \;1\;1\;0\;\cdots$&\\
\cmidrule{1-1}
$\phantom{ .\;0\;}\cdots \; 0\; 1^{v-1} 1\;0\;  1^{u-3} \;0\;1\;1\;\cdots$\phantom{ .}&\hspace{-0.12cm}\mbox{\raisebox{+0.42cm}[0pt]{.}}
\end{tabular}
\end{center}
Note that %$\ell(3k)=\ell+2$ and
\begin{align}\label{eq:4}
L_{v+u+2}(3k) = 0\, 1^{v-1}\, 0^2 \,1^{u-2} \, 0 \, 1.
\end{align}
Thus we get
\begin{align*}
t_{kn} &\equiv s(L_{v+u+2}(3k)) + t_{3k} + v +(u-3) +2+ s(U_{u+1}(k)) + t_k\\
 &\equiv (v-1) + (u-2) + 1 + t_{3k}+ v +(u-3) +2+ u+ t_k\\
 &\equiv \modd{1 + t_{3k} +t_k} {2},
\end{align*}
and we obtain $f(k)<k$.

Now we consider the case $u = 2$.  In order to complete the proof of the
lemma, it remains to show that $f(k)<k$ for integers $k$ with  $U_{3}(k)=1^2 0$ and $L_{v+4}(k) = 0 1^v 0 1^2$ for an even positive integer $v$. If $U_{4}(k)=1^2 0 1$, then we set $n=2^{\ell-4}+1$. Here we get
\[
(kn)_2 = U_{\ell-(v+4)}(k)\, 1\, 0^{v-1} \, 1 \, 0^3 \, L_{\ell-4}(k),
\]
as illustrated below:
\begin{center}
\begin{tabular}{ll}
$\phantom{ .}\cdots \;0 \; 1^{v-1} \;1\;0 \;1\;1$&\\
$\phantom{ .\cdots \;0 \; 1^{v-1}} \;1\;1 \;0\;1\;\cdots$&\\
\cmidrule{1-1}
$\phantom{ .}\cdots \;1 \; 0^{v-1} \;1\;0 \;0\;0\;\cdots$\phantom{ .}&\hspace{-0.12cm}\mbox{\raisebox{+0.42cm}[0pt]{,}}
\end{tabular}
\end{center}
and we obtain
\begin{align*}
t_{kn} &\equiv s(L_{v+4}(k)) + t_{k} + 1+1 + s(U_{4}(k)) + t_k\\
&\equiv (v + 2)  + 3 \equiv \modd{1} {2}.
\end{align*}
If $U_{5}(k)=1^2 0^3$, we set $n=2^{\ell-4}+2^{\ell-5}+1$. It follows that
\[
(kn)_2 = U_{\ell+2-(v+4)}(3k)\, 1\, 0^{v-1} \, 1 \, 0^2 \,1\, L_{\ell-5}(k),
\]
as illustrated below:
\begin{center}
\begin{tabular}{ll}
$\phantom{ .}\cdots \;0 \;1\; 1^{v-2} \;1\;0 \;1\;1$&\\
$\phantom{ .0\;}\cdots \;0 \; 1^{v-2} \;1\;1\;0 \;1\;1$&\\
$\phantom{ .0\;\cdots \;0 \; 1^{v-2}} \;1\;1\;0 \;0\;0\;\cdots$&\\
\cmidrule{1-1}
$\phantom{ .0\;}\cdots \;1 \; 0^{v-2} \;0\;1\;0 \;0\;1\;\cdots$\phantom{ .}&\hspace{-0.12cm}\mbox{\raisebox{+0.42cm}[0pt]{,}}
\end{tabular}
\end{center}
and we obtain
\begin{align*}
t_{kn} &\equiv s(L_{v+4}(3k)) + t_{3k} + 1+1 +1+ s(U_{5}(k)) + t_k\\
&\equiv v+t_{3k} + 1  + 2 + t_k \equiv \modd{1 + t_{3k} + t_k} {2}.
\end{align*}
Here we used Eq.~\eqref{eq:4} and we get $f(k)<k$.

If $U_{5}(k)=1^2 0^2 1$, we set $n=2^{\ell-3}+2^{\ell-1}+1$. It is easy to
see that $\ell(5k)=\ell+2$ or $\ell(5k)=\ell+3$. We have
\[
(kn)_2 = U_{\ell(5k)-(v+4)}(5k)\, 1\, 0^{v+3} \, L_{\ell-5}(k),
\]
as illustrated below:
\begin{center}
\begin{tabular}{ll}
$\phantom{ .}\cdots \;0 \;1\;1\; 1^{v-2} \;0 \;1\;1$&\\
$\phantom{ .0\;0\;}\cdots \;0 \; 1^{v-2} \;1\;1\;0 \;1\;1$&\\
$\phantom{ .0\;0\;\cdots \;0 \; 1^{v-2}} \;1\;1\;0 \;0\;1\;\cdots$&\\
\cmidrule{1-1}
$\phantom{ .0\;0\;}\cdots \;1 \; 0^{v-2} \;0\;0\;0 \;0\;0\;\cdots$\phantom{ .}&\hspace{-0.12cm}\mbox{\raisebox{+0.42cm}[0pt]{.}}
\end{tabular}
\end{center}
Since
\[
L_{v+4}(5k) = 0 \, 1^{v-2} \,0^2 \, 1^3,
\]
we obtain
\begin{align*}
t_{kn} &\equiv s(L_{v+4}(5k)) + t_{5k} + 1+ s(U_{5}(k)) + t_k\\
&\equiv (v+1)+t_{5k} + 1  + 3 + t_k \equiv \modd{1 + t_{5k} + t_k} {2}.
\end{align*}
Again, the considerations of Section~\ref{section:idea} show that $f(k)<k$.
\end{proof}

%%%%%%%%%%%%%%%%%%%%%%%%%%%%%%%%%%%%%%%%%%%%%%%%%%%%%%%%%%%%%%%%%%%%%%%%%%%%%%%%%%%%%%%%%%%%%%%%%%%%%%%%%%%%%%%%%%%%%%%%%%%%%%%%%%
%%%%%%%%%%%%%%%%%%%%%%%%%%%%%%%%%%%%%%%%%%%%%%%%%%%%%%%%%%%%%%%%%%%%%%%%%%%%%%%%%%%%%%%%%%%%%%%%%%%%%%%%%%%%%%%%%%%%%%%%%%%%%%%%%%

\begin{lem}\label{lem:3}
Let $k\in\N$. If there exists an even integer $u\ge2$ and a positive integer $r\ne u$
such that $L_{u+2}(k)=0^21^u$ and $U_{r+1}(k) = 1^r 0$, then we have $f(k)< k$.
\end{lem}

\begin{proof}
Let $\ell=\ell(k)$.
If $r<u$, we set $n= 2^{\ell-r-1}+1<k$. In exactly the same manner
as at the beginning of the proof of Lemma~\ref{lem:2}
(see Eq.~\eqref{eq:3}), we see that
$t_{kn}=1$ and thus, $f(k)<k$.

If $r>u$, we set $n=2^{\ell-u-1}+1<k$. Then we get
\[
(kn)_2 = U_{\ell-(u+2)}(k)\, 1\, 0\, 1^{u-1}\, 0\, L_{\ell-(u+1)}(k),
\]
as illustrated below:
\begin{center}
\begin{tabular}{ll}
$\phantom{ .}\cdots \;0 \;0\;1^{u-1} \;1$&\\
$\phantom{ .\cdots \;0} \;1\;1^{u-1} \;1\;\cdots$&\\
\cmidrule{1-1}
$\phantom{ .}\cdots \;1 \;0\;1^{u-1} \;0\;\cdots$\phantom{ .}&\hspace{-0.12cm}\mbox{\raisebox{+0.42cm}[0pt]{.}}
\end{tabular}
\end{center}
Similarly as before, we have $s(U_{\ell-(u+2)}(k)) \equiv s(L_{u+2}(k))
+ t_k \equiv \modd{u + t_k} {2}$ and $s(L_{\ell-(u+1)}(k)) \equiv
s(U_{u+1})+t_k \equiv u+1+ t_k$ (note that $u+1 \le r$). We obtain
\[
t_{kn} \equiv u + t_k + 1 + (u-1) + (u+1)+ t_k \equiv u+ 1 \equiv \modd{1} {2}.
\]
This shows the desired result.
\end{proof}

%%%%%%%%%%%%%%%%%%%%%%%%%%%%%%%%%%%%%%%%%%%%%%%%%%%%%%%%%%%%%%%%%%%%%%%%%%%%%%%%%%%%%%%%%%%%%%%%%%%%%%%%%%%%%%%%%%%%%%%%%%%%%%%%%%
%%%%%%%%%%%%%%%%%%%%%%%%%%%%%%%%%%%%%%%%%%%%%%%%%%%%%%%%%%%%%%%%%%%%%%%%%%%%%%%%%%%%%%%%%%%%%%%%%%%%%%%%%%%%%%%%%%%%%%%%%%%%%%%%%%
\begin{lem}\label{lem:4}
Let $k\in\N$. If there exists an even integer $u\ge2$ and a positive integer $s<u-1$ such that $L_{u+2}(k)=0^21^u$ and $U_{u+s+1}(k) = 1^u 0^s1$, then we have $f(k)< k$.
\end{lem}

\begin{proof}
Let $\ell=\ell(k)$ and set $n= 2^{\ell-1} + 2^{u-1}+1$. Since $k$ is odd and starts with at least two $1$'s, we see that $n<k$. We have
\[
(n k)_2 = U_{\ell-(u+2)}(k) \, 1 \, 0^{u+s+1} L_{\ell(km)-(s+u+1)}(m k),
\]
where $m=2^{u-1}+1<k$ and $\ell(km)= \ell+u$, as illustrated below:
\begin{center}
\begin{tabular}{ll}
$\phantom{ .}\cdots \;0 \;0\;1^{u-1}\; 1$&\\
$\phantom{ .\cdots \;0 \;0}\;1^{u-1}\; 1\;0^s\;1\;\cdots$&\\
$\phantom{ .\cdots \;0 \;0\;1^{u-1}}\; 1\;1^s\;1\;\cdots\phantom{ .}$&\\
\cmidrule{1-1}
$\phantom{ .}\cdots \;1 \;0\;0^{u-1}\; 0\;0^s\;\cdots$&\hspace{-0.12cm}\mbox{\raisebox{+0.42cm}[0pt]{.}}
\end{tabular}
\end{center}
We have
\[
U_{s+u+1}(km) = 1\, 0^{u-1} 1 0^s.
\]
In particular, we obtain
\[
s(L_{\ell(km)-(s+u+1)}(m k)) \equiv s(U_{s+u+1}(km)) + t_{km}
\equiv \modd{t_{km}} {2},
\]
and we get
\[
t_{n k} \equiv s(L_{u+2}(k))+ t_k + 1+ t_{km} \equiv u + 1+ t_k +
t_{km} \equiv \modd{1+ t_k + t_{km}} {2}.
\]
As before, we get $f(k)<k$, which proves Lemma~\ref{lem:4}.
\end{proof}

%%%%%%%%%%%%%%%%%%%%%%%%%%%%%%%%%%%%%%%%%%%%%%%%%%%%%%%%%%%%%%%%%%%%%%%%%%%%%%%%%%%%%%%%%%%%%%%%%%%%%%%%%%%%%%%%%%%%%%%%%%%%%%%%%%
%%%%%%%%%%%%%%%%%%%%%%%%%%%%%%%%%%%%%%%%%%%%%%%%%%%%%%%%%%%%%%%%%%%%%%%%%%%%%%%%%%%%%%%%%%%%%%%%%%%%%%%%%%%%%%%%%%%%%%%%%%%%%%%%%%

\begin{lem}\label{lem:5}
Let $k\in\N$. If there exists an even integer $u\ge2$ and a positive integer $t\ge 2$ such that $L_{u+t+2}(k)=0 1 0^t 1^u$ and $U_{2u-1}(k) = 1^u 0^{u-1}$, then we have $f(k)< k$.
\end{lem}

\begin{proof}
Let $\ell=\ell(k)$ and let us first assume that $2\le t  \le u-1$.
We set $n=2^{\ell-(u+t)}+1<k$. Then we get
\[
(kn)_2 = U_{\ell-(u+t+2)}(k)\, 1\, 0^{t+1}\, 1^{u-(t+1)}\, 0\, 1^{t} \, L_{\ell-(u+t)}(k),
\]
as illustrated below:
\begin{center}
\begin{tabular}{ll}
$\phantom{ .}\cdots \;0 \;1\;0^{t}\; 1^{u-(t+1)}\; 1 \; 1^t$&\\
$\phantom{ .\cdots \;0 \;1}\;1^{t}\; 1^{u-(t+1)}\; 1 \; 0^t\;\cdots$&\\
\cmidrule{1-1}
$\phantom{ .}\cdots \;1 \;0\;0^{t}\; 1^{u-(t+1)}\; 0 \; 1^t\;\cdots$\phantom{ .}&\hspace{-0.12cm}\mbox{\raisebox{+0.42cm}[0pt]{.}}
\end{tabular}
\end{center}
Since $u< u+t \le 2u-1$, we have
\[
s(L_{\ell-(u+t)}(k)) \equiv s(U_{u+t}(k)) +t_k \equiv \modd{u + t_k} {2}.
\]
This implies
\[
t_{kn} \equiv s(L_{u+t+2}(k)) +t_k + 1+ (u-(t+1)) + t+ u + t_k \equiv u +1  \equiv \modd{1} {2}.
\]
If $t=u$ and $U_{2u}(k) = 1^u 0^{u-1} 1$, we again set
$n=2^{\ell-(u+t)}+1<k$. This time we can write
\[
(kn)_2 = U_{\ell-(u+t+2)}(k)\, 1\, 0^{t+u+1}\, L_{\ell-(u+t)}(k),
\]
as illustrated below:
\begin{center}
\begin{tabular}{ll}
$\phantom{ .}\cdots \;0 \;1\;0^{t}\; 1^{u-1}\; 1 $&\\
$\phantom{ .\cdots \;0 \;1}\;1^{t}\; 0^{u-1}\; 1 \;\cdots$&\\
\cmidrule{1-1}
$\phantom{ .}\cdots \;1 \;0\;0^{t}\; 0^{u-1}\; 0 \;\cdots$\phantom{ .}&\hspace{-0.12cm}\mbox{\raisebox{+0.42cm}[0pt]{,}}
\end{tabular}
\end{center}
and we get
\[
t_{kn} \equiv s(L_{u+t+2}(k))+t_k + 1+ s(U_{u+t}(k)) +t_k \equiv (u +1) +1+ (u+1)  \equiv \modd{1} {2}.
\]
Alternatively, if $t=u$ and $U_{2u+1}(k) = 1^u 0^{u} a$ for some $a\in\{0,1\}$, then we set $n = 2^{\ell-(t+u+1)}+1<k$. Since we have (recall that $\bar{a}=1-a$)
\[
(kn)_2 = U_{\ell-(u+t+2)}(k)\, 1\,0\,  1^{t-1}\,a\, {\bar{a}}^{u}\, L_{\ell-(u+t+1)}(k),
\]
as illustrated below,
\begin{center}
\begin{tabular}{ll}
$\phantom{ .}\cdots \;0 \;1\;0^{t-1}\;0\; 1^{u-1}\; 1 $&\\
$\phantom{ .\cdots \;0} \;1\;1^{t-1}\;0\; 0^{u-1}\; a \;\cdots$&\\
\cmidrule{1-1}
$\phantom{ .}\cdots \;1 \;0\;1^{t-1}\;a\; {\bar{a}}^{u-1}\; \bar{a} \;\cdots$\phantom{ .}&\hspace{-0.12cm}\mbox{\raisebox{+0.42cm}[0pt]{,}}
\end{tabular}
\end{center}
and we finally obtain
\begin{align*}
t_{kn} &\equiv s(L_{u+t+2}(k))+t_k + 1+  (t-1) + a + a u  +   s(U_{u+t+1}(k)) +t_k\\
&\equiv (u +1) +1+ (t-1) + a(u+1) + (u + a)\\
&\equiv \modd{1} {2}.
\end{align*}
This shows the desired result.
\end{proof}

%%%%%%%%%%%%%%%%%%%%%%%%%%%%%%%%%%%%%%%%%%%%%%%%%%%%%%%%%%%%%%%%%%%%%%%%%%%%%%%%%%%%%%%%%%%%%%%%%%%%%%%%%%%%%%%%%%%%%%%%%%%%%%%%%%
%%%%%%%%%%%%%%%%%%%%%%%%%%%%%%%%%%%%%%%%%%%%%%%%%%%%%%%%%%%%%%%%%%%%%%%%%%%%%%%%%%%%%%%%%%%%%%%%%%%%%%%%%%%%%%%%%%%%%%%%%%%%%%%%%%

\begin{lem}\label{lem:6}
Let $k\in\N$. If there exists an even integer $u\ge2$ and positive integers $t\ge 2$ such that $L_{u+t+2}(k)=11 0^t 1^u$ and $U_{2u-1}(k) = 1^u 0^{u-1}$, then we have $f(k)< k$.
\end{lem}

\begin{proof}
Let $\ell=\ell(k)$. First, we consider the case $2\le t \le u-1$.
Set $n=2^{\ell-(t+u-1)} + 2^{\ell-(t+u)} + 1<k$. Then we have
\[
(kn)_2 = U_{\ell+2-(t+u+2)}(3k)\, 1\, 0^{t-1}  \, 1 \, 0 \, 1^{u-(t+1)} \, 0\,  1^{t-2} \, 0 \, 1 \, L_{\ell-(u+t)}(k),
\]
as illustrated below:
\begin{center}
\begin{tabular}{ll}
$\phantom{ .}\cdots \;1 \;1\;0\; 0^{t-2}\;0\;1\; 1^{u-(t+1)}\;1\; 1^{t-2}\; 1$&\\
$\phantom{ .\;1}\cdots \;1 \;1\; 0^{t-2}\;0\;0\; 1^{u-(t+1)}\;1\; 1^{t-2}\; 1\;1$&\\
$\phantom{ .\;1\cdots \;1 \;1}\; 1^{t-2}\;1\;1\; 1^{u-(t+1)}\;1\; 0^{t-2}\; 0\;0\;\cdots$&\\
\cmidrule{1-1}
$\phantom{ .\;1}\cdots \;1 \;0\; 0^{t-2}\;1\;0\; 1^{u-(t+1)}\;0\; 1^{t-2}\; 0\;1\;\cdots$\phantom{ .}&\hspace{-0.12cm}\mbox{\raisebox{+0.42cm}[0pt]{.}}
\end{tabular}
\end{center}
Noting that $\ell(3k)=\ell+2$, we obtain
\begin{align*}
t_{kn} &\equiv s(L_{u+t+2}(3k)) + t_{3k} + (u-(t+1)) + (t-2) + s(U_{u+t}(k)) + t_k +3\\
 &\equiv s(L_{u+t+2}(3k)) + t_{3k} +u + u +t_k \\
 &\equiv \modd{s(L_{u+t+2}(3k)) + t_{3k} +t_k} {2}.
\end{align*}
Since
\begin{align}\label{eq:2}
L_{u+t+2}(3k) = 0 \, 1 \, 0^{t-2} \,1 \, 0 \, 1^{u-2} \, 0 \, 1,
\end{align}
we have $t_{kn} \equiv \modd{1 + t_{3k}+t_k} {2}$ and consequently $f(k)<k$.

If $t=u$, then we have to consider three different cases. If $U_{2u}(k) = 1^u 0^{u-1} 1$ we again set $n = 2^{\ell-(t+u-1)}+ 2^{\ell-(t+u)}+1<k$.
This time we get
\[
(kn)_2 = U_{\ell+2-(t+u+2)}(3k)\, 1\, 0^{t}  \, 1^u \, 0\, L_{\ell-(2u)}(k),
\]
as illustrated below,
\begin{center}
\begin{tabular}{ll}
$\phantom{ .}\cdots \;1 \;1\;0\; 0^{t-1}\;1\; 1^{u-1}$&\\
$\phantom{ .\;1}\cdots \;1 \;1\; 0^{t-1}\;0\; 1^{u-1}\;1$&\\
$\phantom{ .\;1\cdots \;1 \;1}\; 1^{t-1}\;1\; 0^{u-1}\;1\;\cdots$&\\
\cmidrule{1-1}
$\phantom{ .\;1}\cdots \;1 \;0\; 0^{t-1}\;1\; 1^{u-1}\;0\;\cdots$\phantom{ .}&\hspace{-0.12cm}\mbox{\raisebox{+0.42cm}[0pt]{,}}
\end{tabular}
\end{center}
which yields
\begin{align*}
t_{kn} &\equiv s(L_{u+t+2}(3k)) + t_{3k} + 1 + u + s(U_{2u}(k)) + t_k\\
 &\equiv s(L_{u+t+2}(3k)) + t_{3k} + 1 + u +(u+1) +t_k \\
 &\equiv \modd{s(L_{u+t+2}(3k)) + t_{3k} +t_k} {2} .
\end{align*}
As above (see Eq.~\eqref{eq:2}), we get that $f(k)<k$.
If $t=u$ and $U_{2u+1}(k) = 1^u 0^{u} 1$, then we set $n= 2^{\ell-(t+u)} +  2^{\ell-(t+u)-1} +1<k$. We have
\[
(kn)_2 = U_{\ell+2-(t+u+2)}(3k)\, 1^2 \, 0^{t}  \, 1^{u-1} \, 0 \,L_{\ell-(2u+1)}(k),
\]
as illustrated below.
\begin{center}
\begin{tabular}{ll}
$\phantom{ .}\cdots \;1 \;1\;0\; 0^{t-1}\;1\; 1^{u-1}$&\\
$\phantom{ .\;1}\cdots \;1 \;1\; 0^{t-1}\;0\; 1^{u-1}\;1$&\\
$\phantom{ .\;1\cdots \;1} \;1\; 1^{t-1}\;0\; 0^{u-1}\;1\;\cdots$&\\
\cmidrule{1-1}
$\phantom{ .\;1}\cdots \;1 \;1\; 0^{t-1}\;0\; 1^{u-1}\;0\;\cdots$\phantom{ .}&\hspace{-0.12cm}\mbox{\raisebox{+0.42cm}[0pt]{.}}
\end{tabular}
\end{center}
Note that the dots in first and second line of the figure have to be erased if $k=51$. (The binary representation of $51$ is given by $(51)_2 = 110011$.)
We get
\begin{align*}
t_{kn} &\equiv s(L_{u+t+2}(3k)) + t_{3k} + 2 + (u-1) + s(U_{2u+1}(k)) + t_k\\
 &\equiv s(L_{u+t+2}(3k)) + t_{3k} +2+ (u-1) + (u+1) + t_k \\
 &\equiv \modd{s(L_{u+t+2}(3k)) + t_{3k} +t_k} {2} .
\end{align*}
Using Eq.~\eqref{eq:2}, we obtain $f(k)<k$.
If $t=u$ but $U_{2u+1}(k) = 1^u 0^{u+1}$, then we choose $n=2^{\ell-1} + 2^{u}+1$. Since $k$ is odd and starts with at least two $1$'s, we again obtain that $n<k$. This leads us to
\[
(kn)_2 = U_{\ell-(u+2)}(k)\, 1 \, 0 \,  1^{u-1} \,0 \, 1^{u-1} \,L_{\ell+u-(2u)}(km),
\]
where $m=2^{u}+1<k$, as illustrated below:
\begin{center}
\begin{tabular}{ll}
$\phantom{ .}\cdots \;0 \;0\;1^{u-1}\;1$&\\
$\phantom{ .\cdots \;0} \;1\;1^{u-1}\;0\; 0^{u-1}\; 0\;\cdots$&\\
$\phantom{ .\cdots \;0 \;1\;1^{u-1}}\;1\; 1^{u-1}\; 0\;\cdots$\phantom{ .}\\
\cmidrule{1-1}
$\phantom{ .}\cdots \;1 \;0\;1^{u-1}\;0\; 1^{u-1}\;\cdots$&\hspace{-0.12cm}\mbox{\raisebox{+0.42cm}[0pt]{.}}
\end{tabular}
\end{center}
Note that $\ell(km)=\ell+u$ and
\[
U_{2u}(km) = 1^{2u}.
\]
We obtain
\begin{align*}
t_{kn} &\equiv s(L_{u+2}(k)) + t_{k} + 1 + (u-1) + (u-1) + s(U_{2u}(km)) + t_{km}\\
 &\equiv u + t_{k} +1 +2u+t_{km} \\
 &\equiv \modd{1 + t_{k} +t_{km}} {2},
\end{align*}
which shows $f(k)<k$ in this case, too.

In order to prove the lemma, it remains to consider the case $t>u$. We set $n= 2^{\ell-1} + 2^{\ell-u-1}+1<k$ and we get
\[
(kn)_2 = U_{\ell+u-(2u+1)}(km)\, 1 \, 0^{u-1} \,  1^{u-1} \,0 \, 1 \,L_{\ell-(u+1)}(k)),
\]
as illustrated below,
\begin{center}
\begin{tabular}{ll}
$\phantom{ .}\cdots \;0 \; 1^{u-1}\;1$&\\
$\phantom{ .}\cdots \;0 \; 0^{u-1}\;0\; 1^{u-2}\; 1\; 1\;$&\\
$\phantom{ .\cdots \;0 \; 0^{u-1}}\;1\; 1^{u-2}\; 1\; 0\;\cdots$&\\
\cmidrule{1-1}
$\phantom{ .}\cdots \;1 \; 0^{u-1}\;1\; 1^{u-2}\; 0\; 1\;\cdots$\phantom{ .}&\hspace{-0.12cm}\mbox{\raisebox{+0.42cm}[0pt]{,}}
\end{tabular}
\end{center}
where $m=2^{u}+1<k$. Since
\[
L_{2u+1}(km) = 0\, 1^{2u},
\]
we obtain
\begin{align*}
t_{kn} &\equiv s(L_{2u+1}(km)) + t_{km} + 1 + (u-1) + 1 + s(U_{u+1}(k)) + t_{k}\\
 &\equiv 2u + t_{km} + u-1 +u +t_{k} \\
 &\equiv \modd{1 + t_{k} +t_{km}} {2}.
\end{align*}
The same argument as before finally shows the desired result.
\end{proof}

\section{Proof of the Theorem~\ref{thm:1}}\label{section:theorem}

\begin{proof}[Proof of Theorem~\ref{thm:1}]
As already noted in Section~\ref{section:idea}, we have $f(k)=f(2k)$ for all $k\ge1$. Consequently, it suffices to show that $f(k)\le k+4$ for odd integers $k$.

If $k=2^{2r+1}-1$, $r\ge0$, then $(k)_2 = 1^{2r+1}$ and we trivially have $f(k)=1$. If $k=2^{2r}-1$, $r\ge1$, then we will show that $f(k)=k+4$.
In order to do this, we note that  the binary sum-of-digits function $s_2$ satisfies the relation
\[
s_2(a 2^k - b) = s_2(a-1) + k - s_2(b-1)
\]
for all positive integers $a,b,k$ with $1\le b < 2^k$.
Thus we have for all $1\le m \le k$,
\[
t_{km} \equiv s_2(2^{2r} m  - m) = s_2(m-1) + 2r - s_2(m-1) \equiv
\modd{0} {2}.
\]
If $m=2^{2r}$ or $m=2^{2r}+2$, then we clearly have $t_{km}=0$ since $t_{km/2}=0$. If $m=2^{2m}+1$, then $km = 2^{4r}-1$ and consequently $t_{km}=0$. If $m=2^{2r}+3$, then
\[
t_{km} \equiv s_2(2^{4r} + 2^{2r+1}-3) \equiv 1 + (2r+1) - s_2(2) \equiv
\modd{1} {2},
\]
which finally proves $f(k)=k+4$ for $k=2^{2r}-1$, $r\ge1$. If $k$ is a positive integer different from ones already considered, then there exist positive integers $r,s,t$ and $u$ such that
\[
U_{r+s+1}(k) = 1^r\, 0^s\, 1 \qquad \mbox{ and } \qquad L_{t+u+1} = 1 \, 0^t \, 1^u.
\]
If $u$ is odd, then Lemma~\ref{lem:1} implies that $f(k)\le k$ where equality occurs if and only if $k= 2^{r}+1$ for some $r\ge2$. If $u$ is even but $t=1$ or $r\ne u$, then Lemma~\ref{lem:2} and Lemma~\ref{lem:3} imply $f(k)<k$. Let us assume that $u$ is even, $t\ge2$ and $r=u$. Then there exists $a\in\{0,1\}$, such that
\[
L_{t+u+2} = a\, 1 \, 0^t \, 1^u.
\]
If $s<u-1$, then Lemma~\ref{lem:4} implies $f(k)<k$. Contrarily, if $s\ge u-1$, then Lemma~\ref{lem:5} (if $a=0$) or Lemma~\ref{lem:6} (if $a=1$)  yield $f(k)<k$. Hence we have for all positive integers $k$,
\[
f(k) \le k+4,
\]
where equality occurs if and only if $k=2^{2r}-1$, $r\ge1$. Note that an even positive integer $2m$ cannot satisfy $f(2m)=2m+4$, since we then would get $f(m)=f(2m) = 2m +4\le  m+4$ and consequently $m\le0$. Moreover, we see that for odd integers $k$ there exist no solutions to the equation
\[
f(k)=k+\alpha
\]
for $\alpha= 0,1,2,3$,
except in the case $\alpha=0$ where
we have $f(k) = k$ if and only if $k=2^r+1$ for some $r\ge2$ or $k=1$.
If $k=2m$ is even, then $f(2m) = 2m+\alpha$ implies $f(m) = 2m +\alpha
\le m+4$. Hence this can only happen if $m\le 4-\alpha$. We see that
there exist no solutions to $f(k)=k+\alpha$ for $\alpha = 2$ and $
\alpha=3$, there are no even solutions for $\alpha=0$ and the only
solution to $f(k)=k+1$ is $k=6$ (compare with Eq.~\eqref{eq:5}). This
finally proves Theorem~\ref{thm:1}.
\end{proof}

\begin{rem}
  By a similar case analysis it might be possible to prove that $$\min\{n: t_{kn}=0\}\leq k+2.$$ However, it does not seem
  possible to obtain this bound in a direct way from the bound~(\ref{thm:1}).
\end{rem}

\section{Some weak general results}\label{section:general}

Given the generality of Gelfond's theorem, it is natural to try to bound the minimal $n$ such that $n \equiv \modd{a} {k}$ and $s_b(n) \equiv \modd{c} {r}$.
Here we only get a weaker upper bound.

\begin{prop}\label{prope}
Let $b, r, k$ be positive integers with $\gcd(b-1, r) = 1$,
and let $c$ be any integer.  Then there exists
a non-negative integer $n < b^r k$ such that
$s_b (kn) \equiv \modd{c} {r}$.
\label{easy}
\end{prop}

\begin{proof}
We claim that if $1 \leq k \leq b^t$, then $s_b (k (b^t - 1)) = (b-1)t$.
To see this, note that
$$s_b(kb^t -k)=s_b(k-1)+(b-1)t-s_b(k-1)=(b-1)t.$$
Let $s$ be the smallest integer such that $k \leq b^s$.
Then $b^{s-1} < k$.
Choose $t \in \lbrace s, s+1, \ldots, s+r-1\}$
such that $(b-1) t \equiv \modd{c} {r}$.   This is possible
since $\gcd(b-1,r) = 1$.
Then $s_b( k (b^t - 1)) = (b-1) t \equiv \modd{c} {r}$,
as desired.  Furthermore, $b^t \leq b^{s+r-1} \leq b^r b^{s-1}
< b^r k$.  Thus we can take $n = b^t - 1$.
\end{proof}

\begin{cor}
Let $b, r, k$ be positive integers with $\gcd(b-1, r) = 1$,
and let $a, c$ be any integers.  Then there exists an integer
$n < b^{r+1} k^3$ such that
$n \equiv \modd{a} {k}$ and $s_b (n) \equiv \modd{c} {r}$.
\end{cor}

\begin{proof}
Without loss of generality we can assume $0 \leq a < k$.
As in the proof of Proposition~\ref{easy} let $s$ be the smallest
integer such that $b^s \geq k$, so $b^{s-1} < k$.
From Proposition~\ref{easy} we know that there exists an
integer $t$ such that $s_b (k (b^t - 1)) \equiv \modd{(c-a)} {r}$,
and $b^t < b^r k$.  Then clearly
$s_b ( k b^s (b^t - 1) + a) \equiv \modd{c} {r}$, so we can take
$n = k b^s (b^t-1) + a$.   Then $n < b^{r+1} k^3$.
\end{proof}

In the setting of Proposition~\ref{prope} we conjecture that a
similar phenomenon takes place as we have seen in the case of
the classical Thue-Morse sequence.

\begin{conj}
Let $b, r$ be positive integers with $\gcd(b-1, r) = 1$, and let
$c$ be any integer.
There exists a constant $C$, depending only on $b$ and $r$ such that
for all $k \geq 1$ there exists $n \leq k + C$ with
$s_b(kn) \equiv \modd{c} {r}$.  Furthermore, we can take
$C \leq b^{r+c}$.
\end{conj}

\section{Acknowledgments}

We are very grateful to Jorge Buescu for asking the initial question
that led us to the results in this paper.  We also thank Jean-Paul Allouche
for allowing us to reproduce his remarks here.

\bibliographystyle{elsarticle-num}

\begin{thebibliography}{9}

\bibitem{Allouche:1982}
J-.P. Allouche.
\newblock Somme des chiffres et transcendance.
\newblock {\it Bull. Soc. Math. France} {\bf 110} (1982), 279--285.

\bibitem{DartygeLucaStanica:2009}
C. Dartyge, F. Luca, P. St{\u{a}}nic{\u{a}}.
\newblock On digit sums of multiples of an integer.
\newblock {\it J. Number Theory} {\bf 129} (2009), 2820--2830.

\bibitem{Fine:1965}
N. J. Fine.
\newblock The distribution of the sum of digits $({\rm mod}~p)$.
\newblock {\it Bull. Amer. Math. Soc.} {\bf 71} (1965), 651--652.

\bibitem{Gelfond:1968}
A. O. Gelfond.
\newblock Sur les nombres qui ont des propri\'et\'es additives et multiplicatives donn\'ees.
\newblock {\it Acta Arithmetica} {\bf 13} (1968), 259--265.
\newblock Electronically available at {\tt http://matwbn.icm.edu.pl/ksiazki/aa/aa13/aa13115.pdf}.

\bibitem{Newman:1969}
D. J. Newman.
\newblock On the number of binary digits in a multiple of three.
\newblock {\it Proc. Amer. Math. Soc.} {\bf 21} (1969), 719--721.

\end{thebibliography}

\end{document}